\documentclass{amsart}
\usepackage{graphicx}
\theoremstyle{definition}

\theoremstyle{remark}

\numberwithin{equation}{section}

\begin{document}
\title[Characterizations of all-derivable points in nest
algebras]{Characterizations of all-derivable points in nest
algebras{\footnote{This work is supported by the National Natural
Science Foundation of China (No 10771191)}}}%
\author{Jun Zhu and Sha Zhao}%
\address{Institute of Mathematics, Hangzhou Dianzi
University, Hangzhou 310018, PR China}%
\email{junzhu@yahoo.cn}%

\thanks{}%
\begin{abstract}
Let $\mathcal{A}$ be an operator algebra on a Hilbert space. We say
that an element $G\in {\mathcal{A}}$ is an all-derivable point of
${\mathcal{A}}$ if every derivable linear mapping $\varphi$ at $G$
(i.e. $\varphi(ST)=\varphi(S)T+S\varphi(T)$ for any $S,T\in
alg{\mathcal{N}}$ with $ST=G$) is a derivation. Suppose that
$\mathcal{N}$ is a nontrivial complete nest on a Hilbert space $H$.
We show in this paper that $G\in {alg\mathcal{N}}$ is an
all-derivable point if and only if $G\neq0$.
\end{abstract}
\maketitle
\section{1. Introduction and preliminaries}
Let $\mathcal{H}$ and $\mathcal{K}$ be two Hilbert spaces. $B(H,K)$
stands for the set of all bounded linear operators from $H$ into
$K$, and is abbreviated to $B(H,H)$ to $B(H)$. Let $\mathcal{A}$ be
an operator subalgebra in $B(H)$. We say that a linear mapping
$\varphi$ from $\mathcal{A}$ into itself is a derivable mapping at
$G$ if $\varphi(ST)=\varphi(S)T+S\varphi(T)$ for any $S,T\in
{\mathcal{A}}$ with $ST=G$. We say that an element $G\in
\mathcal{A}$ is an all-derivable point of $\mathcal{A}$ if every
derivable linear mapping $\varphi$ at $G$ is a derivation. Let $N$
be a closed subspace in $H$. We use the symbols $I_{N}$ to denote
the unit operator in $B(N)$. If $\mathcal{N}$ is a complete nest on
$H$, then the nest algebra $alg{\mathcal{N}}$ is the set of all
operators which leave every member of $\mathcal{N}$ invariant.

As well known, derivations are very important maps both in theory
and applications. In general there are two directions in the study
of the local actions of derivations of operator algebras. One is the
well known local derivation problem (see [9]). The other is to study
conditions under which derivations of operator algebras can be
completely determined by the action on some sets of operators. It is
obvious that a liner map is a derivation if and only if it is
derivable at all points. It is natural and interesting to ask the
question whether or not a linear map is a derivation if it is
derivable only at one given point.

We describe some of the results related to ours.  Jing and Lu [7]
showed that every derivable mapping $\varphi$ at $0$ with
$\varphi(I)=0$ on nest algebras is a derivation. Hou and Qi [6] got
that every idempotent with the range in $\mathcal{N}$ is an
all-derivations in ${alg\mathcal{N}}$, where $\mathcal{N}$ is a
complete nest on Banach space. Li, Pan and Xu [11] proved that every
derivable mapping $\varphi$ at $0$ with $\varphi(I)=0$ on CSL
algebras is a derivation. Zhu and Xiong in [16] showed that every
element $G\in {\mathcal{TM}}_{n}$ is an all-derivable point of
${\mathcal{TM}}_{n}$ if and only if $G\neq 0$, where
${\mathcal{TM}}_{n}$ is the algebra of all $n\times n$ upper
triangular matrices. For other relative reference, see
[1-5,8-10,12-15].

It is the aim of this paper to prove that an operator $G\in
{alg\mathcal{N}}$ is an all-derivable point of the nest algebra
$alg\mathcal{N}$ if and only if $G\neq 0$, where $\mathcal{N}$ be a
nontrivial complete nest on a Hilbert space $H$ (i.e. there exists
an element $N\in {\mathcal{N}}$ with $N\neq \{0\}, H$).
\section{2. All-derivable points in nest algebras}
Suppose that $N$ is a closed subspace in $H$. In this section, we
always use the symbols $I_{N}$ to denote the
unit operator on $N$. We easily prove the following two lemmas.\\
 {\bf{Lemma 2.1}} {\it Let
$\mathcal{N}$ be a complete nest on a Hilbert space $H$. Then $G$ is
an all-derivation point of $alg{\mathcal{N}}$ if and only if
$\lambda G$ is of $alg{\mathcal{N}}$ for any real number
$\lambda\neq0$.\\}  {\bf Proof.} $\varphi$ is a derivable mapping at
$G$ $\Leftrightarrow$ $\varphi(ST)=\varphi(S)T+S\varphi(T)$ for any
$S,T\in {\mathcal{A}}$ with $ST=G$ $\Leftrightarrow$
$\varphi(\lambda ST)=\varphi(\lambda S)T+\lambda S\varphi(T)$ for
any $S,T\in {\mathcal{A}}$ with $ST= G$ $\Leftrightarrow$ $\varphi(
ST)=\varphi(S)T+S\varphi(T)$ for any $S,T\in {\mathcal{A}}$ with
$ST= \lambda G$.$\Box$\\
{\bf{Lemma 2.2}} {\it Let $\mathcal{A}$ be an operator subalgebra
with unit operator $I$ in $B(H)$, and let $\varphi$ is a linear
mapping from $\mathcal{A}$ into itself. Suppose that $\varphi(X)=0$
for any invertible operator $X\in \mathcal{A}$, Then $\varphi\equiv
0$.}\\ {\bf Proof.} For arbitrary $X\in alg{\mathcal{N}}$, there
exists a real number $\lambda>\parallel X\parallel$. Then both
$\lambda I-X$ and $2\lambda I-X$ are two invertible operators. So
$\varphi(\lambda I-X)=0$ and $\varphi(2\lambda I-X)=0$. It follows
from the linearity of $\varphi$ that $\varphi(X)=0$.$\Box$

The following is our main theorem in this paper.\\
 {\bf{Theorem 2.3}} {\it
Let $\mathcal{N}$ be a nontrivial complete nest on a Hilbert space
$H$. Then an operator $G\in {alg\mathcal{N}}$ is an all-derivable
point if and only if $G\neq 0$.}
\\
{\bf Proof.} Suppose that $\varphi$ is a derivable linear mapping at
$G\neq0$ from $alg{\mathcal{N}}$ into itself. We only need to prove
that $\varphi$ is a derivation. Let $N\in {\mathcal{N}}$ with
$\{0\}\subset N\subset H$. Then all $2\times 2$ operator matrices
always are represented as relative to the orthogonal decomposition
$H=N\oplus N^{\perp}$ in the proof of this theorem. Thus we may
write $$ G=\left[\begin{array}{cc}
D & E \\
0 & F
\end{array}\right]$$
where $D\in alg{\mathcal{N}}_{N}$, $E\in
alg{\mathcal{N}}_{N^{\perp}}$ and $F\in B(N^{\perp},N)$ (
${\mathcal{N}}_{N}=\{M\cap N: M\in {\mathcal{N}}\}$ and
${\mathcal{N}}_{N}=\{M\cap N^{\perp}: M\in {\mathcal{N}}\}$).
Without loss of generality, we may assume that $\parallel
D\parallel<1$ by Lemma 2.1. For arbitrary $X\in
alg{\mathcal{N}}_{N}$, $Y\in B(N^{\perp},N)$ and $Z\in
alg{\mathcal{N}}_{N^{\perp}}$, we write
$$\left\{\begin{array}{ccc}\varphi(\left[\begin{array}{cc}
X & 0 \\
0 & 0
\end{array}\right])=\left[\begin{array}{cc}
A_{11}(X) & A_{12}(X) \\
0 & A_{22}(X)
\end{array}\right],\\~\\
\varphi(\left[\begin{array}{cc}
0 & Y \\
0 & 0
\end{array}\right])=\left[\begin{array}{cc}
B_{11}(Y) & B_{12}(Y) \\
0 & B_{22}(Y)
\end{array}\right],\\~\\
\varphi(\left[\begin{array}{cc}
0 & 0 \\
0 & Z
\end{array}\right])=\left[\begin{array}{cc}
C_{11}(Z) & C_{12}(Z) \\
0 & C_{22}(Z)
\end{array}\right].\end{array}\right.$$
Obviously, $A_{ij},B_{ij}$ and $C_{ij}$ $(i,j=1,2, i\leq j)$ are
linear mappings on $alg{\mathcal{N}}_{N}$, $B(N^{\perp},N)$ and
$alg{\mathcal{N}}_{N^{\perp}}$, respectively.

Step 1. We show that $A_{11}(\cdot)$ is a derivable mapping at $D$.
For arbitrary $X,U\in alg{\mathcal{N}}$ with $XU=D$, taking
$S=\left[\begin{array}{cc}
\lambda^{-1}X & \lambda E \\
0 & \lambda F
\end{array}\right], T=\left[\begin{array}{cc}
\lambda U & 0 \\
0 & \lambda^{-1}I_{N^{\perp}}
\end{array}\right]\in alg{\mathcal{N}}$ for any real number $\lambda>0$, then $ST=G$. It follows
that
\begin{eqnarray}
&&\left[\begin{array}{cc} A_{11}(D)+B_{11}(E)+C_{11}(F) &
A_{12}(D)+B_{12}(E)+C_{12}(F) \\ 0 & A_{22}(D)+B_{22}(E)+C_{22}(F)
\end{array}\right]\\\nonumber&=&\varphi(G)=\varphi(S)T+S\varphi(T)\\\nonumber&=&\left[\begin{array}{cc}
\lambda^{-1}A_{11}(X)+\lambda B_{11}(E)+\lambda C_{11}(F) & * \\
0 & *
\end{array}\right]\left[\begin{array}{cc}
\lambda U & 0 \\
0 & \lambda^{-1}I_{N^{\perp}}
\end{array}\right]\\\nonumber
\end{eqnarray}\begin{eqnarray}\nonumber&=&\left[\begin{array}{cc}
\lambda^{-1}X & \lambda E \\
0 & \lambda F
\end{array}\right]\left[\begin{array}{cc}
\lambda A_{11}(U)+\lambda^{-1}C_{11}(I) & * \\ 0 & *
\end{array}\right]\\\nonumber
&=&\left[\begin{array}{cc}
A_{11}(X)U+\lambda^{2}B_{11}(E)U+\lambda^{2}C_{11}(F)U & * \\
+XA_{11}(U)+\lambda^{-2}XC_{11}(I)&~\\ 0 & *
\end{array}\right].
\end{eqnarray}
It follows from the matrix equation that
$$A_{11}(D)+B_{11}(E)+C_{11}(F)=A_{11}(X)U+\lambda^{2}B_{11}(E)U+\lambda^{2}C_{11}(F)U+XA_{11}(U)+\lambda^{-2}XC_{11}(I)$$
The above equation implies that
$$A_{11}(D)+B_{11}(E)+C_{11}(F)=A_{11}(X)U+XA_{11}(U);$$
$$B_{11}(E)U+C_{11}(F)U=0;XC_{11}(I)=0.$$
Furthermore $C_{11}(I)=0$ and $B_{11}(E)+C_{11}(F)=0$. Thus
$A_{11}(D)=A_{11}(X)U+XA_{11}(U)$, i.e. $A_{11}(\cdot)$ is a
derivable mapping at $D$.

Step 2. We shows that $C_{11}(W)=0$ and $B_{11}(V)=0$ for any $W\in
alg{\mathcal{N}}_{N^{\perp}}$ and $V\in B(N^{\perp},N)$.

Letting $S=\left[\begin{array}{cc}
X & Y \\
0 & Z
\end{array}\right], T=\left[\begin{array}{cc}
U & V \\
0 & W
\end{array}\right]\in alg{\mathcal{N}}$ with $ST=G$, then
$XU=D$, $XV+YW=E$ and $ZW=F$. Since $\varphi$ is a derivable mapping
at $G$ on $alg\mathcal{N}$, we have
\begin{eqnarray}
&&\left[\begin{array}{cc} A_{11}(D)+B_{11}(E)+C_{11}(F) &
A_{12}(D)+B_{12}(E)+C_{12}(F) \\\nonumber 0 &
A_{22}(D)+B_{22}(E)+C_{22}(F)
\end{array}\right]\\&=&\varphi(G)=\varphi(S)T+S\varphi(T)\\\nonumber&=&\left[\begin{array}{cc}
A_{11}(X)+B_{11}(Y) & A_{12}(X)+B_{12}(Y)  \\
+C_{11}(Z)&+C_{12}(Z)\\
 0 & A_{22}(X)+B_{22}(Y)+C_{22}(Z)
\end{array}\right]\left[\begin{array}{cc}
U & V \\
~&~\\ 0 & W
\end{array}\right]\\\nonumber
&=&\left[\begin{array}{cc}
X & Y \\
~&~\\ 0 & Z
\end{array}\right]\left[\begin{array}{cc}
A_{11}(U)+B_{11}(V) & A_{12}(U)+B_{12}(V) \\
+C_{11}(W)&+C_{12}(W)\\
0 & A_{22}(U)+B_{22}(V)+C_{22}(W)
\end{array}\right].
\end{eqnarray}
The above equation implies the following three equations hold.
\begin{eqnarray}
&&\nonumber
A_{11}(D)+B_{11}(E)+C_{11}(F)\\&=&A_{11}(X)U+B_{11}(Y)U+C_{11}(Z)U\\\nonumber&&+XA_{11}(U)+XB_{11}(V)+XC_{11}(W);
\end{eqnarray}
\begin{eqnarray}
&&A_{12}(D)+B_{12}(E)+C_{12}(F)\nonumber\\&=&A_{11}(X)V+B_{11}(Y)V+C_{11}(Z)V\\\nonumber&&+A_{12}(X)W+B_{12}(Y)W+C_{12}(Z)W\\\nonumber
&&+XA_{12}(U)+XB_{12}(V)+XC_{12}(W)\\\nonumber&&+YA_{22}(U)+YB_{22}(V)+YC_{22}(W);
\end{eqnarray}
\begin{eqnarray}
&&\nonumber
A_{22}(D)+B_{22}(E)+C_{22}(F)\\&=&A_{22}(X)W+B_{22}(Y)W+C_{22}(Z)W\\\nonumber&&+ZA_{22}(U)+ZB_{22}(V)+ZC_{22}(W);
\end{eqnarray}
Note that $A_{11}$ is a derivable mapping at $D$. So we have
$A_{11}(XU)=A_{11}(D)=A_{11}(X)U+XA_{11}(U)$ for any $X,U\in
alg{\mathcal{N}}_{N}$. By Eq. (2.3), we have
\begin{eqnarray}
 B_{11}(E)+C_{11}(F)=B_{11}(Y)U+C_{11}(Z)U+XB_{11}(V)+XC_{11}(W).
\end{eqnarray}
Taking $X=D, U=I_{N}, Y=E, V=0, Z=F$ and $W=I_{N^{\perp}}$ in Eq.
(2.6), then $DC_{11}(I_{N^{\perp}})=0$. Letting $X=D, U=I_{N}, Y=-E,
V=0, Z=-F$ and $W=-I_{N^{\perp}}$ in Eq. (2.6), then
$B_{11}(E)+C_{11}(F)=0$. It follows that
$$0=B_{11}(Y)U+C_{11}(Z)U+XB_{11}(V)+XC_{11}(W).$$
For arbitrary $V\in B(N^{\perp},N)$, if we put $X=I_{N}, U=D,
W=I_{N^{\perp}}, Z=F$ and $Y=E-V$ in the above equation, then
$B_{11}(V)(I-D)=0$. Since $\parallel D\parallel<1$ implies that
$I-D$ is invertible, $B_{11}(V)=0$. For arbitrary $W\in
alg{\mathcal{N}}_{N^{\perp}}$, if we put $X=I_{N}, Y=0, Z=FW^{-1},
U=D$ and $V=E$ in the above equation, then
$$C_{11}(FW^{-1})D+C_{11}(W)=0.$$ Replacing $W$ by $2W$ in the above equation, we have
$$\frac{1}{2}C_{11}(FW^{-1})D+2C_{11}(W)=0.$$
The above two equations implies that $C_{11}(W)=0$ for any
invertible operator $W\in alg{\mathcal{N}}_{N^{\perp}}$. By Lemma
2.2, $C_{11}(W)\equiv0$.

Step 3. We show that $A_{22}(X)=0$ and $B_{22}(Y)=0$ for any $X\in
alg{\mathcal{N}}_{N}$ and $Y\in B(N^{\perp},N)$.

Letting $X=D$, $U=I_{N}$, $Y=E$, $V=0$, $Z=F$ and
$W=I_{{N}^{\perp}}$ in Eq. (2.5), then
\begin{eqnarray}F(A_{22}(I_{N})+C_{22}(I_{{N}^{\perp}}))=0.\end{eqnarray} On the other hand, for
arbitrary $Y\in B(N^{\perp},N)$ and real number $\lambda\neq 0$,
putting $X=\lambda I_{N}$, $U=\lambda^{-1}D$, $V=Y+\lambda^{-1}E$,
$Z=-\lambda^{-1}F$ and $W=-\lambda I_{{N}^{\perp}}$ in Eq. (2.5),
then
\begin{eqnarray}
A_{22}(D)+B_{22}(E)\nonumber&=&-\lambda^{2}A_{22}(I_{N})-\lambda
B_{22}(Y)-\lambda^{-2}FA_{22}(D)\\\nonumber&-&\lambda^{-1}FB_{22}(Y)-\lambda^{-2}FB_{22}(E)+FC_{22}(I_{{N}^{\perp}}).
\end{eqnarray}
The above equation implies that
\begin{eqnarray}
&& A_{22}(D)+B_{22}(E)=FC_{22}(I_{{N}^{\perp}}); A_{22}(I_{N})=0;\\
&&FA_{22}(D)+FB_{22}(E)=0; B_{22}(Y)=0.
\end{eqnarray}
It follows from Eq. (2.7)-(2.9) that $FC_{22}(I_{{N}^{\perp}}))=0$
and
$$A_{22}(D)=FC_{22}(I_{{N}^{\perp}})=0.$$ Bring Eq. (2.9) and the above
equation to Eq. (2.5) to get
$$C_{22}(F)=A_{22}(X)W+C_{22}(Z)W+ZA_{22}(U)+ZC_{22}(W).$$
For any invertible operator $X\in alg{\mathcal{N}}_{N}$, taking
$U=X^{-1}D$, $V=0$, $W=I_{N^{\perp}}$, $Y=E$ and $Z=F$ in the above
equation, then we have
$$A_{22}(X)+FA_{22}(X^{-1}D)=0.$$
Replacing $X$ by $2X$, we have
$$2A_{22}(X)+\frac{1}{2}FA_{22}(X^{-1}D)=0.$$
The above two equations implies that $A_{22}(X)=0$ for any
invertible operator $X\in alg{\mathcal{N}}_{N}$. By Lemma 2.2,
$A_{22}(X)\equiv0$.

Step 4. We show that $A_{12}(X)=-XC_{12}(I_{{N}^{\perp}})$ and
$C_{12}(W)=C_{12}(I_{{N}^{\perp}})W$ for any $X\in
alg{\mathcal{N}}_{N}$ and $W\in alg{\mathcal{N}}_{N^{\perp}}$. For
arbitrary invertible operator $X\in alg{\mathcal{N}}_{N}$, taking
$U=X^{-1}D$, $V=0$, $Y=\lambda E$, $Z=\lambda F$ and
$W=\lambda^{-1}I_{{N}^{\perp}}$ in Eq. (2.4), then
\begin{eqnarray}
A_{12}(D)&=&\lambda^{-1}A_{12}(X)
+XA_{12}(X^{-1}D)\\\nonumber&&+\lambda^{-1}XC_{12}(I_{{N}^{\perp}})+\lambda
EA_{22}(X^{-1}D)+EC_{22}(I_{{N}^{\perp}});
\end{eqnarray}
Eq. (2.10) implies that
$$A_{12}(X)=-XC_{12}(I_{{N}^{\perp}}).$$
It follows from Lemma 2.2 that $A_{12}(X)=-XC_{12}(I_{{N}^{\perp}})$
for any $X\in alg{\mathcal{N}}$.

For arbitrary invertible operator $W\in
alg{\mathcal{N}}_{{N}^{\perp}}$, taking $X=I_{N}$, $U=D$, $Y=0$,
$V=E$ and $Z=FW^{-1}$ in Eq. (2.4), then
\begin{eqnarray}
C_{12}(F)=A_{11}(I_{N})E+A_{12}(I_{N})W+C_{12}(FW^{-1})W+C_{12}(W).
\end{eqnarray}
Replacing $W$ by $2W$ in the above equation, we have
\begin{eqnarray}
C_{12}(F)=A_{11}(I_{N})E+2A_{12}(I_{N})W+C_{12}(FW^{-1})W+2C_{12}(W).
\end{eqnarray}
Combining Eq. (2.11) with Eq. (2.12), we obtain
$$A_{12}(I_{N})W+C_{12}(W)=0.$$
Letting $W=I_{{N}^{\perp}}$ in the above equation, then
$A_{12}(I_{N})=-C_{12}(I_{{N}^{\perp}})$. Hence
$C_{12}(W)=C_{12}(I_{{N}^{\perp}})W$. It follows from Lemma 2.2 that
$C_{12}(W)=C_{12}(I_{{N}^{\perp}})W$ for any $W\in
alg{}\mathcal{{N}^{\perp}}$.

Step 5. We show that both $A_{11}(I_{N})=0$ and
$C_{22}(I_{{N}^{\perp}})=0$.

For arbitrary $V\in B(N^{-1},N)$, taking $X=I_{N}$, $U=D$, $Y=E-V$,
$Z=F$ and $W=I_{{N}^{\perp}}$ in Eq. (2.4), by the results of Step
2-3, we have
\begin{eqnarray}
A_{11}(I_{N})V=VC_{22}(I_{{N}^{\perp}})-EC_{22}(I_{{N}^{\perp}}).
\end{eqnarray}
Letting $V=0$ in Eq. (2.13), then $EC_{22}(I_{{N}^{\perp}})=0$.
Hence $ A_{11}(I_{N})V=VC_{22}(I_{{N}^{\perp}})$ for any $V\in
B(N^{\perp},N)$. For arbitrary $x\in N$ and $y\in N^{-1}$, we have
$x\otimes y\in B(N^{\perp},N)$. So $ A_{11}(I_{N})x\otimes
y=x\otimes C_{22}(I_{{N}^{\perp}})^{*}y$. Thus there exists a
complex number $\alpha$ such that $A_{11}(I_{N})=\alpha I_{N}$ and
$C_{22}(I_{{N}^{\perp}})=\alpha I_{{N}^{\perp}}$. Note that
$FC_{22}(I_{{N}^{\perp}})=0$ and $EC_{22}(I_{{N}^{\perp}})=0$. So
$\alpha F=0$ and $\alpha E=0$. On the other hand, taking $X=I_{N}$,
$U=D$, $V=E$, $Y=0$, $Z=I$ and $W=F$ in Eq. (2.3), then
$A_{11}(I_{N})D=0$, i.e. $\alpha D=0$. Since $G\neq 0$, $\alpha=0$.
It follows that $A_{11}(I_{N})=0$ and $C_{22}(I_{{N}^{\perp}})=0$.

Step 6. We show that both $A_{11}(\cdot)$ and $C_{22}(\cdot)$ are
derivations.

For arbitrary $Y\in B(N^{\perp},N)$ and invertible operator $W\in
alg{\mathcal{N}}_{N^{\perp}}$, taking $X=I_{N}$, $U=D$, $V=E-YW$ and
$Z=FW^{-1}$ in Eq. (2.4), then
\begin{eqnarray}B_{12}(YW)=B_{12}(Y)W+YC_{22}(W).\end{eqnarray}
For any $W_{1},W_{2}\in alg{\mathcal{N}}_{N^{\perp}}$, we have
\begin{eqnarray}B_{12}(YW_{1}W_{2})=B_{12}(Y)W_{1}W_{2}+YC_{22}(W_{1}W_{2}).\end{eqnarray}
On the other hand,
\begin{eqnarray}&&B_{12}(YW_{1}W_{2})=B_{12}(YW_{1})W_{2}+YW_{1}C_{22}(W_{2})\\\nonumber
&=&B_{12}(Y)W_{1}W_{2}+YC_{22}(W_{1})W_{2}+YW_{1}C_{22}(W_{2})
.\end{eqnarray} The Eq. (2.15) and Eq. (2.16) implies that
$$C_{22}(W_{1}W_{2})=C_{22}(W_{1})W_{2}+W_{1}C_{22}(W_{2}).$$
Hence $C_{22}(\cdot)$ is a derivation.

For arbitrary $V\in B(N^{\perp},N)$ and invertible operator $X\in
alg{\mathcal{N}_{N}}$, taking $X^{-1}D$, $Y=E-XV$,
$W=I_{{N}^{\perp}}$ and $Z=F$, then
\begin{eqnarray}B_{12}(XV)=A_{11}(X)V+XB_{12}(V).\end{eqnarray}
Similarly, we can prove that $A_{11}(\cdot)$ is a derivation.

Step 7. We show that $\varphi$ is a derivation. For arbitrary
$S=\left [\begin{array}{ccc}
  X & Y \\
  0 & Z \\
\end{array}\right ]$ and $T=\left
[\begin{array}{ccc}
  U & V \\
  0 & W \\
\end{array}\right ]$ in
$alg{\mathcal{N}}$,  we only need to prove that
$\varphi(ST)=\varphi(S)T+S\varphi(T)$. By Eqs. (2.14), (2.17) and
the results of Step 2-6, we easily calculate
\begin{eqnarray}
\nonumber&&\varphi(S)T+S\varphi(T)\\\nonumber&=&\left
[\begin{array}{ccc}
  A_{11}(X) & A_{12}(X)+B_{12}(Y)+C_{12}(Z) \\
  0 & C_{22}(Z) \\
\end{array}\right ]\left [\begin{array}{ccc}
  U & V \\
  0 & W \\
\end{array}\right ]\\\nonumber
&&+\left [\begin{array}{ccc}
  X & Y \\
  0 & Z \\
\end{array}\right ]\left
[\begin{array}{ccc}
  A_{11}(U) & A_{12}(U)+B_{12}(V)+C_{12}(W) \\
  0 & C_{22}(W) \\
\end{array}\right ]\\\nonumber
&=&\left [\begin{array}{ccc}
  A_{11}(X)U+XA_{11}(U) & A_{11}(X)V+A_{12}(X)W+B_{12}(Y)W \\
  ~& +C_{12}(Z)W+XA_{12}(U)+XB_{12}(V)\\
  ~&+XC_{12}(W)+YC_{22}(W)\\
  ~&~\\
  0 & C_{22}(Z)W+ZC_{22}(W) \\
\end{array}\right ]\\\nonumber
&=&\left [\begin{array}{ccc}
  A_{11}(XU) & A_{12}(X)W+B_{12}(YW)+C_{12}(Z)W \\
  ~& +XA_{12}(U)+B_{12}(XV)+XC_{12}(W)\\
  ~&~\\
  0 & C_{22}(ZW) \\
\end{array}\right ]\\\nonumber
&=&\left [\begin{array}{ccc}
  A_{11}(XU) & A_{12}(X)W+B_{12}(YW+XV) \\
  ~& +C_{12}(Z)W+XA_{12}(U)+XC_{12}(W)\\
  ~&~\\
  0 & C_{22}(ZW) \\
\end{array}\right ]\\\nonumber
&=&\left [\begin{array}{ccc}
  A_{11}(XU) & -XC_{12}(I_{{N}^{\perp}})W+B_{12}(YW+XV) \\
  ~& +C_{12}(I_{{N}^{\perp}})ZW-XUC_{12}(I_{{N}^{\perp}})+XC_{12}(I_{{N}^{\perp}})W\\
  ~&~\\
  0 & C_{22}(ZW) \\
\end{array}\right ]\\\nonumber
&=&\left [\begin{array}{ccc}
  A_{11}(XU) & B_{12}(YW+XV)+C_{12}(I_{{N}^{\perp}})ZW-XUC_{12}(I_{{N}^{\perp}}) \\
  ~&~\\
  0 & C_{22}(ZW) \\
\end{array}\right ]\\\nonumber
&=&\left [\begin{array}{ccc}
  A_{11}(XU) & B_{12}(YW+XV)+C_{12}(ZW)+A_{12}(XU) \\
  ~&~\\
  0 & C_{22}(ZW) \\
\end{array}\right ]\\\nonumber&=&\varphi(ST).
\end{eqnarray}
This completes the proof.
$\Box$\\
\bibliographystyle{amsplain}\bibliography{Referenes}
{\bf{Referenes}}\\ 1. M. Bre$\check{s}$ar, Characterization of
derivations on some normed algebras with involution, J.Algebra 152
(1992), 454-462. MR 1194314 (94e:46098) \\
2. M. Bre$\check{s}$ar, P. $\check{S}$emrl,  Mappings which preserve
idempotents, local automorphisms, and local
  derivations. Canad. J. Math. 45 (1993), 483-496. MR 1222512 (94k:47054)\\
3. R.L. Crist, Local derivations on operator algebras, J. Funct.
Anal. 135 (1996), 76-92.\\
4. K.R. Davidson, Nest algebras, Research Notes in Math. No. 191,
Longman Sci. \& Tech., Wiley \& Sons, New York, 1988. MR 0972978 (90f:47062) \\
5. J.A. Erdos, Operator of finite rank in nest algebras, J. London
Math. Soc. 43 (1968), 391-397.  MR 0230156 (37 $\#$5721) \\
6. J.C. Hou, X.F. Qi, Characterizations of derivations of
  Banach space nest algebras: all-derivable point, Linear Algebra
  Appl. 432 (2010), 3183-3200.\\
7. W. Jing, S.J. Lu, P.T. Li,  Characterizations of
  derivations on some operator algebras, Bull. Austral. Math. Soc.
  66 (2002), 227-232. MR 1932346 (2003f:47059)\\
8. R.V. Kadison, Local derivations, J. Algebra 130 (1990),
  494-509. MR 1051316 (91f:46092)\\
9. D.R. Larson, A.R. Sourour, Local derivations and local
automorphisms of ${\mathcal{B}}(X)$,
  Operator Algebras and Applications, Proc. Symp. Pure Math. 51
  (1990), 187-194. MR 1077437 (91k:47106) \\
10.  D.R. Larson, Nest algebras and similarity
  transformations, Annals of Mathematics 121 (1985), 409-427. MR0794368 (86j:47061)\\
11.  J.K. Li, Z.D. Pan, H. Xu,  Characterizations of isomorphisms
and derivations of some algebras. J. Math. Anal. Appl.
  332 (2007), 1314-1322. MR 2324339 (2008e:47093)\\
12. P. $\check{S}$emrl, Local automorphisms and derivations on
$B(H)$, Proc. Amer. Math. Soc.
  125 (1997), 2677-2680. MR 1415338 (98e:46082)\\
13. J. Zhu, All-derivable points of operator algebras. Linear
Algebra Appl.  427 (2007), 1-5. MR 2353150 (2008g:47141) \\
14.  J. Zhu, C.P. Xiong, Derivable mappings at unit
  operator on
  nest algebras. Linear Algebra Appl. 422 (2007),
  721-735. MR 2305152 (2008a:47123)\\
15. J. Zhu, C.P. Xiong, Bilocal derivations of standard operator
algebras, Proc. Amer. Math. Soc.
  125 (1997), 1367-1370. MR 1363442 (97g:46091)\\
16. J. Zhu, C.P. Xiong, R.Y. Zhang, All-derivable points in the
algebra of all upper triangular matrices
  Linear Algebra Appl. 429(4) (2008), 804-818. MR2428131 (2009m:47200) \\

\end{document}